\documentclass[11pt]{article}
\usepackage{amssymb,amsmath,latexsym}
\usepackage{pdfsync}

\begin{document}

\begin{doublespace}

\def\1{{\bf 1}}

\def\sA {{\cal A}} \def\sB {{\cal B}} \def\sC {{\cal C}}
\def\sD {{\cal D}} \def\sE {{\cal E}} \def\sF {{\cal F}}
\def\sG {{\cal G}} \def\sH {{\cal H}} \def\sI {{\cal I}}
\def\sJ {{\cal J}} \def\sK {{\cal K}} \def\sL {{\cal L}}
\def\sM {{\cal M}} \def\sN {{\cal N}} \def\sO {{\cal O}}
\def\sP {{\cal P}} \def\sQ {{\cal Q}} \def\sR {{\cal R}}
\def\sS {{\cal S}} \def\sT {{\cal T}} \def\sU {{\cal U}}
\def\sV {{\cal V}} \def\sW {{\cal W}} \def\sX {{\cal X}}
\def\sY {{\cal Y}} \def\sZ {{\cal Z}}

\def\bA {{\mathbb A}} \def\bB {{\mathbb B}} \def\bC {{\mathbb C}}
\def\bD {{\mathbb D}} \def\bE {{\mathbb E}} \def\bF {{\mathbb F}}
\def\bG {{\mathbb G}} \def\bH {{\mathbb H}} \def\bI {{\mathbb I}}
\def\bJ {{\mathbb J}} \def\bK {{\mathbb K}} \def\bL {{\mathbb L}}
\def\bM {{\mathbb M}} \def\bN {{\mathbb N}} \def\bO {{\mathbb O}}
\def\bP {{\mathbb P}} \def\bQ {{\mathbb Q}} \def\bR {{\mathbb R}}
\def\bS {{\mathbb S}} \def\bT {{\mathbb T}} \def\bU {{\mathbb U}}
\def\bV {{\mathbb V}} \def\bW {{\mathbb W}} \def\bX {{\mathbb X}}
\def\bY {{\mathbb Y}} \def\bZ {{\mathbb Z}}
\def\R {{\mathbb R}} \def\RR {{\mathbb R}}
\def\n{{\bf n}}

\newcommand{\expr}[1]{\left( #1 \right)}
\newcommand{\cl}[1]{\overline{#1}}
\newtheorem{thm}{Theorem}[section]
\newtheorem{lemma}[thm]{Lemma}
\newtheorem{defn}[thm]{Definition}
\newtheorem{prop}[thm]{Proposition}
\newtheorem{corollary}[thm]{Corollary}
\newtheorem{remark}[thm]{Remark}
\newtheorem{example}[thm]{Example}
\numberwithin{equation}{section}
\def\ee{\varepsilon}
\def\qed{{\hfill $\Box$ \bigskip}}
\def\NN{{\cal N}}
\def\AA{{\cal A}}
\def\MM{{\cal M}}
\def\BB{{\cal B}}
\def\CC{{\cal C}}
\def\LL{{\cal L}}
\def\DD{{\cal D}}
\def\FF{{\cal F}}
\def\EE{{\cal E}}
\def\QQ{{\cal Q}}
\def\RR{{\mathbb R}}
\def\R{{\mathbb R}}
\def\L{{\bf L}}
\def\K{{\bf K}}
\def\S{{\bf S}}
\def\A{{\bf A}}
\def\E{{\mathbb E}}
\def\F{{\bf F}}
\def\P{{\mathbb P}}
\def\N{{\mathbb N}}
\def\eps{\varepsilon}
\def\wh{\widehat}
\def\wt{\widetilde}
\def\pf{\noindent{\bf Proof.} }
\def\beq{\begin{equation}}
\def\eeq{\end{equation}}
\def\bee{\begin{equation}}
\def\eee{\end{equation}}

\title{\Large \bf
Uniform Boundary Harnack Principle for
Rotationally Symmetric L\'evy processes in General Open Sets}
\author{{\bf Panki Kim}\thanks{This research was supported by Basic Science Research Program
through the National Research Foundation of Korea(NRF) funded by the Ministry of Education, Science and Technology(0409-20110087).} \quad {\bf Renming Song}\thanks{Research supported in part by a grant from the Simons
Foundation (208236).}
\quad and \quad {\bf Zoran Vondra\v{c}ek}\thanks{Supported in part by the MZOS
grant 037-0372790-2801.}  }

\date{ }
\maketitle

\begin{abstract}
In this paper we prove the uniform boundary Harnack principle in general open sets for harmonic functions
with respect to a large class of rotationally symmetric purely discontinuous L\'evy processes.
\end{abstract}



\section{Introduction}
The boundary Harnack principle for classical harmonic functions is a very deep result in potential
theory and has many important applications in probability theory and analysis.

In the late nineties Bogdan \cite{B} established the boundary Harnack principle
for harmonic functions of rotationally symmetric $\alpha$-stable processes, $\alpha\in (0, 2)$,
in Lipschitz domains. This was the first time that the boundary Harnack principle was established for harmonic
functions with respect to non-local operators (or, equivalently, discontinuous Markov processes).
Since then the result has been generalized in various directions. In \cite{SW} Song and Wu extended the
boundary Harnack principle to harmonic functions with respect to rotationally
symmetric  stable processes in $\kappa$-fat open sets, with the constant depending on the local geometry
near the boundary. The definitive result in the case of rotationally symmetric stable processes was obtained in
\cite{BKK} by Bogdan, Kulczycki and Kwa\'snicki who established the boundary Harnack principle in arbitrary
opens sets with the constant not depending on the open set itself. This type of result is known
as the uniform boundary Harnack principle.
Note that the uniform boundary Harnack principle is not true for Brownian motion.

In another direction, the boundary Harnack principle has been generalized to different classes of
discontinuous processes. In \cite{KSV} the boundary Harnack principle was established for harmonic functions
with respect to a wide class of purely discontinuous subordinate Brownian motions in $\kappa$-fat open sets,
with an extension obtained in \cite{KSV3}.
In \cite{KSV4} (see also, \cite{CKSV, KSV2}) the boundary Harnack inequality was established for
harmonic functions of subordinate Brownian motions with Gaussian components.

The purpose of this paper is to generalize the main results from
\cite{BKK, KSV, KSV3} and prove the uniform
boundary Harnack principle for harmonic functions with respect to a large class of
rotationally symmetric purely discontinuous L\'evy processes in arbitrary open sets.
The class of processes treated in this paper is larger than the class of processes
treated in \cite{KSV, KSV3}.
 The processes considered in this paper need not be subordinate
Brownian motions. Even when restricted to subordinate Brownian motions, the assumptions
on the subordinate Brownian motions in this paper are slightly weaker than those in
\cite{KSV, KSV3}.

To be more precise, let $S=(S_t:\, t\ge 0)$ be a subordinator with Laplace exponent $\phi$. We assume that
$\phi$ is a complete Bernstein function  satisfying the following
{\it upper and lower scaling conditions} (see \cite{Z}):

\medskip
\noindent
{\bf (H):}
There exist constants $ \delta_1, \delta_2 \in (0,1)$, $a_1, a_2>0$ and $R_0>0$ such that
\begin{eqnarray*}
{\rm (LSC)} \qquad &
\phi(\lambda r) \ge a_1\lambda^{\delta_1} \phi(r), \quad &\lambda \ge 1, r \ge 1/R_0^2
\\
{\rm (USC)} \qquad &
\phi(\lambda r) \le a_2 \lambda^{\delta_2} \phi(r), \quad &\lambda \ge 1, r \ge 1/R_0^2.
\end{eqnarray*}
Note that it follows from {\rm (USC)} that $\phi$ has no drift.

Let $W=(W_t:\, t\ge 0)$ be a Brownian motion in $\R^d$, $d\ge 1$, independent of the subordinator $S$.
The subordinate Brownian motion $Y=(Y_t:\, t\ge 0)$ is defined by $Y_t:=
W_{S_t}$.
The L\'evy measure of the process $Y$ has a density given by
$J(x)=j(|x|)$ where
\begin{equation}\label{e:ld4s}
j(r):=\int^{\infty}_0(4\pi t)^{-d/2}e^{-r^2/(4t)}\mu(t)dt, \qquad r>0
\end{equation}
and $\mu(t)$ is the L\'evy density of $S$. Note that the function
$r\mapsto j(r)$ is continuous and decreasing
on $(0, \infty)$.

We will assume that $X$ is a purely discontinuous rotationally
symmetric L\'evy process with L\'evy exponent
$\Psi(\xi)$.  Because of rotational symmetry, the function $\Psi$ depends on $|\xi|$ only, and by a slight abuse of notation we write $\Psi(\xi)=\Psi(|\xi|)$. We further assume that the L\'evy measure of $X$ has a density $J_X$. Then
$$
\E_x\left[e^{i\xi\cdot(X_t-X_0)}\right]=e^{-t\Psi(|\xi|)},
\quad \quad \mbox{ for every } x\in \R^d \mbox{ and } \xi\in \R^d,
$$
with
\begin{equation}\label{e:psi}
\Psi(|\xi|)= \int_{\R^d}(1-\cos(\xi\cdot y))J_X(y)dy.
\end{equation}
We assume that
$J_X$ is continuous on $\R^d\setminus\{0\}$ and that
there is a constant $\gamma>1$ such that
\begin{equation}\label{e:psi1}
\gamma^{-1} j(|y|)\le J_X(y) \le \gamma j(|y|), \quad \mbox{for all } y\in \R^d\, .
\end{equation}
Clearly \eqref{e:psi1} implies that
\begin{equation}\label{e:psi2}
\gamma^{-1} \phi(|\xi|^2)\le \Psi(|\xi|) \le \gamma \phi(|\xi|^2), \quad \mbox{for all } \xi\in \R^d\, .
\end{equation}
For a Greenian open set $D\subset \R^d$, we will use $K_D$ to denote the Poisson kernel of $X$ in $D \times \overline{D}^c$  (see \eqref{PK} below).
The goal of this paper is to establish the following result:

\begin{thm}\label{UBHP}
Let $X$ be a purely discontinuous rotationally symmetric L\'evy
process with a continuous L\'evy density $J_X$ satisfying \eqref{e:psi1} where the
complete Bernstein function $\phi$ satisfies {\bf (H)}.
There exists a constant $c= c(\phi, \gamma, d)>0$ such that
\begin{itemize}
    \item[(i)] For every $z_0 \in \R^d$, every open set $D\subset \R^d$, every $r\in (0,1)$
    and for any nonnegative functions $u, v$ in $\R^d$ which are regular harmonic in $D\cap
    B(z_0, r)$ with respect to $X$ and vanish in $D^c \cap B(z_0, r)$, we have
        $$
        \frac{u(x)}{v(x)}\,\le c\,\frac{u(y)}{v(y)}
        $$
        for all $x, y\in D\cap B(z_0, r/2)$.
    \item[(ii)] For every $z_0 \in \R^d$, every
  Greenian open set
     $D\subset \R^d$, every $r\in (0,1)$, we have
    $$
    K_D(x_1, y_1)K_D(x_2, y_2) \le c K_D(x_1, y_2)K_D(x_2, y_1)
    $$
    for all $x_1, x_2 \in D \cap B(z_0, r/2)$ and all $y_1, y_2 \in
    \overline{D}^c \cap B(z_0, r)^c$.
\end{itemize}
\end{thm}

The proof of the above theorem uses some results developed in \cite{KSV3} and several ideas from \cite{BKK}.
In the next section we recall some necessary definitions and results from \cite{KSV3}.
In Section 3 we prove several results about one-dimensional symmetric L\'evy processes that will be needed in
the proof of Theorem \ref{UBHP}.
In Section 4, we present some estimates on the Poisson kernel $K_D$ that are essential for the proof
of Theorem \ref{UBHP}.
The proof of Theorem \ref{UBHP} is given in Section 5, where we also give
an approximate factorization of the Poisson kernel, see  Corollary \ref{c:approx-factor}.
In the last section we relate the assumption {\bf (H)} with the class $OR$ of $O$-regularly
varying functions and sketch the construction of an example of a complete Bernstein function
which satisfies {\bf (H)} but not the assumptions in \cite{KSV3}.

At the meeting ``Foundations of Stochastic Analysis'' held in Banff from
September 18 to 23, 2011, M.~Kwa\'snicki announced that, in a forthcoming
joint paper with K.~Bogdan and T.~Kumagai, they have obtained a version of
the boundary Harnack principle for Hunt processes in metric measure spaces
under rather general conditions.

In this paper we always assume $d \ge 1$. We use the following convention:
The value of the constant $C$ will remain the same throughout this paper,
while $c, c_1, c_2, \cdots$ stand for constants whose values are unimportant and which
may change from location to location.
The dependence of the lower case
constants on the dimension $d$ will not be mentioned explicitly.
The labeling of the constants $c_1, c_2, \cdots$ starts anew in the proof of each result.
The notation $f(t)\asymp g(t)$, $t\to 0$ (respectively $f(t)\asymp g(t)$, $t\to \infty$) means that the quotient $f(t)/g(t)$
stays bounded between two positive constants as $t\to 0$ (respectively $t\to \infty$).

\section{Preliminaries}

Suppose that $S=(S_t: t\ge 0)$ is a subordinator with Laplace exponent $\phi$, that is, $S$ is
a nonnegative L\'evy process with $S_0=0$ and
$$
\E\left[e^{-\lambda S_t}\right]=e^{-t\phi(\lambda)}, \qquad \forall \ t, \lambda>0.
$$
The function $\phi$ can be written in the form
\begin{equation}\label{e:LK}
\phi(\lambda)=b\lambda + \int^\infty_0(1-e^{-\lambda t})\mu(dt)
\end{equation}
where $b\ge 0$ and $\mu$ is a measure on $(0, \infty)$ satisfying $\int^\infty_0
(1\wedge t)\mu(dt)<\infty$. $b$ is called the drift of the subordinator and $\mu$
the L\'evy measure of the subordinator.
The function $\phi$ is a Bernstein function, i.e., it is $C^{\infty}$ and $(-1)^n D^n \phi\ge 0$ for all $n\ge 0$.

Note that, by using \eqref{e:LK} and the elementary inequality $1-e^{-t y}\le
t(1-e^{-y})$ valid for all $t \ge 1$ and all $y>0$, we see that
the Bernstein function $\phi$ satisfies
\begin{equation}\label{e:Berall}
\phi(t\lambda)\le
\lambda\phi(t)
\qquad \text{ for all }
\lambda \ge 1, t >0.
\end{equation}

In this paper we will always assume that $\phi$ is a complete Bernstein function, that is, the L\'evy measure
$\mu$ of $S$ has a
completely monotone density $\mu(t)$, i.e., $(-1)^n D^n\mu\ge 0$ for every non-negative integer $n$.
For basic results on complete
Bernstein functions, we refer our readers to \cite{SSV}.
It follows from \cite[Lemma 2.1]{KSV3} that there exists $c>1$ such that
\begin{equation}\label{e:muatinfty}
\mu(t)\le c \mu(t+1), \qquad t>1.
\end{equation}

The next result will be used to obtain the asymptotic behavior of $\mu(t)$ near the origin.

\begin{prop}\label{p:zahle}
{\rm(\cite[Theorem 7]{Z})}
Suppose that $w$ is a completely monotone function given by
$w(\lambda)=\int^\infty_0e^{-\lambda t} f(t)\, dt,$
where $f$ is a strictly positive decreasing function. Then
$$
f(t)\le \left(1-e^{-1}\right)^{-1} t^{-1}w(t^{-1}), \quad t>0.
$$
If, furthermore, there exist $\delta\in (0, 1)$ and $a, t_0>0$ such
that
\begin{equation}\label{e:zahle}
w(\lambda t)\le a \lambda^{-\delta} w(t), \quad \lambda \ge 1, t\ge 1/t_0,
\end{equation}
then there exists $c=c(f,a,t_0, \delta)>0$ such that
$$
f(t)\ge c t^{-1}w(t^{-1}), \quad t\le t_0.
$$
\end{prop}

From now on we will always assume that the Laplace exponent $\phi$ of $S$
is a complete Bernstein function satisfying {\bf (H)}.

\begin{thm}\label{t:behofmu}
For every $M>0$, there exists $c=c(M, \phi)>1$ such that the L\'evy
density $\mu$ of $S$ satisfies
\begin{equation}\label{e:behofmu}
c^{-1}t^{-1}\phi(t^{-1}) \le  \mu(t)\le c\, t^{-1}\phi(t^{-1})
\quad \text{and} \quad
c^{-1}\phi(t^{-1}) \le  \mu(t, \infty)\le c\, \phi(t^{-1}),  \quad
\forall t \le M\,
 \end{equation}
where $\mu(t,\infty)=\int_t^{\infty}\mu(s)\, ds$ is the tail of the L\'evy measure $\mu$.
\end{thm}

\pf
Let $w(\lambda):=\lambda^{-1}\phi(\lambda)=\int_0^{\infty} e^{-\lambda t}\mu(t,\infty)\, dt$.
The upper scaling condition (USC) implies that $w$ satisfies \eqref{e:zahle} with $\delta=1-\delta_2$
and $t_0=R_0^2$. Hence by Proposition \ref{p:zahle}, there exists a constant $c_1>1$ such that
$$
c_1^{-1}t^{-1}w(t^{-1})\le \mu(t,\infty) \le c_1 t^{-1}w(t^{-1})\, ,\qquad t\le R_0^2\, ,
$$
which immediately implies
\begin{equation}\label{e:v-psi}
c_1^{-1} \phi(t^{-1}) \le \mu(t,\infty) \le c_1 \phi(t^{-1})\, ,\quad t  \le R_0^2\, .
\end{equation}

We proceed to prove the first inequality.
Since $\mu(t/2,\infty)\ge \int_{t/2}^{t}\mu(s)\, ds \ge (t/2)\mu(t)$
by \eqref{e:Berall} and \eqref{e:v-psi}, for
all
$t \in (0, R_0^2]$,
$$
\mu(t)\le 2 t^{-1}\mu(t/2,\infty)\le
2 c_1 t^{-1}\phi((t/2)^{-1})
\le 4c_1 t^{-1}\phi(t^{-1}).
$$
Using {\rm (LSC)} we get that for every $\lambda \ge 1$
\begin{equation}\label{e:nndn}
\phi(s^{-1})=\phi(\lambda(\lambda s)^{-1})\ge
a_1\lambda^{{\delta_1}}\phi((\lambda s)^{-1})\, ,\quad s\le
\frac{R_0^2}{\lambda}\, .
\end{equation}
Fix $\lambda_1:=2^{1/{\delta_1}}((c_1^2 a_1^{-1})\vee 1)^{1/{\delta_1}}\ge1$.
Then, by
\eqref{e:v-psi} and \eqref{e:nndn},  for $s\le (R_0^2\wedge 1)/{\lambda_1}$,
$$
\mu({\lambda_1} s,\infty)\le c_1\phi(({\lambda_1} s)^{-1})\le c_1 a_1^{-1}
{\lambda_1}^{-{\delta_1}}\phi(s^{-1})\le
c_1^2 a_1^{-1}{\lambda_1}^{-{\delta_1}} \mu(s, \infty) \le \frac12 \mu(s,\infty)
$$
by our choice of ${\lambda_1}$. Further,
$$
({\lambda_1}-1)s\mu(s)\ge \int_s^{{\lambda_1} s}\mu(t)\, dt=\mu(s,\infty)-\mu({\lambda_1} s,\infty)\ge \mu(s,\infty)-\frac12 \mu(s,\infty)=\frac12 \mu(s,\infty)\, .
$$
This implies that for all  $t\le (R_0^2\wedge 1)/{\lambda_1}$
$$
\mu(t)\ge \frac{1}{2({\lambda_1}-1)}\, t^{-1} \mu(t,\infty)
\ge \frac{1}{2c_1({\lambda_1}-1)}\, t^{-1} \phi(t^{-1}).
$$
The case $(R_0^2\wedge 1)/{\lambda_1} \le t \le M$ is clear
since the functions we consider are
all positive and continuous on $(0, \infty) $. The proof is now
complete. \qed

A consequence of \eqref{e:behofmu} and (USC) is that for any $K>0$ there exists $c=c(K)>1$ such that
\begin{equation}\label{e:muat0}
\mu(t)\le c \mu(2t), \qquad t\in (0, K).
\end{equation}

Suppose that $W=(W_t:t\ge 0)$ is a
Brownian motion in $\R^d$ with
$$
\E\left[e^{i\xi\cdot(W_t-W_0)}\right]= e^{-t|\xi|^2}, \qquad
\forall\, \xi \in \R^d, t>0\, ,
$$
and that $W$ is independent of $S$. The process $Y=(Y_t:\, t\ge 0)$
defined by $Y_t=W_{S_t}$ is called a subordinate Brownian motion.
It is a rotationally symmetric L\'evy process with
characteristic exponent
$\Phi_Y(\xi)=\phi(|\xi|^2)$, $\xi \in \RR^d$.
Recall that
the L\'evy measure of $Y$ has a density
$J(x)=j(|x|)$ with $j$ given by
\eqref{e:ld4s} and that
$r\mapsto j(r)$ is continuous and decreasing
on $(0, \infty)$.

The following theorem establishes the asymptotic behavior of $j$ near the origin.

\begin{thm}\label{t:Jorigin} It holds that
\begin{equation}\label{e:j-asymptotics}
j(|x|)\asymp \frac{\phi(|x|^{-2})}{|x|^d}
\qquad |x|\to 0.
\end{equation}
\end{thm}
\pf By \eqref{e:Berall},
\begin{equation}\label{e:uv}
\frac{\phi(v)}{v}\le \frac{\phi(u)}{u}\, ,\quad 0<u\le v\, .
\end{equation}
To obtain the upper bound in \eqref{e:j-asymptotics} we write
$$
j(r)=\int_0^{r^2}(4\pi t)^{-d/2} e^{-r^2/(4t)}\mu(t)\, dt +
\int_{r^2}^{\infty}(4\pi t)^{-d/2} e^{-r^2/(4t)}\mu(t)\,
dt:=J_1+J_2\, .
$$
For $r \le 1$, by using \eqref{e:behofmu} in the first inequality and
\eqref{e:uv} in the second, we have
\begin{align*}
&J_1\le  c_1 \int_0^{r^2} (4\pi t)^{-d/2} e^{-r^2/(4t)} t^{-1}
\phi(t^{-1})\, dt \le c_1 \int_0^{r^2} (4\pi t)^{-d/2} e^{-r^2/(4t)}
t^{-1} t^{-1}\frac{\phi(r^{-2})}{r^{-2}}\, dt\\
&\le c_2 r^2 \phi(r^{-2})\int_0^{\infty}t^{-d/2-2}e^{-r^2/(4t)}\,
dt = c_3 r^2 \phi(r^{-2})r^{-d-2}=c_3 r^{-d}\phi(r^{-2})\, .
\end{align*}
Next,
\begin{align*}
&J_2\le  c_4 \int_{r^2}^{\infty}t^{-d/2}\mu(t)\, dt= c_4
\int_{r^2}^{\infty}\left(\frac{d}{2}\int_t^{\infty}
s^{-d/2-1}\, ds\right)\mu(t)\, dt\\
&= c_5 \int_{r^2}^{\infty}\left(\int_{r^2}^s \mu(t)\,
dt\right) s^{-d/2-1}\, ds \le  c_5  \mu(r^2, \infty)
\int_{r^2}^{\infty}s^{-d/2-1}\, ds\le  c_6 r^{-d}\phi(r^{-2})
\end{align*}
where the last inequality  follows from Theorem \ref{t:behofmu}. The
last two displays show that $j(r)\le c_7 r^{-d}\phi(r^{-2})$, for
$r$ small. To prove the converse inequality, we also use Theorem
\ref{t:behofmu} and get that  for $r \le 1$,
\begin{align*}
&j(r)\ge  \int_0^{1}(4\pi s)^{-d/2} e^{-r^2/(4s)}\mu(s)\,
ds=(4\pi)^{-d/2}\int_0^{1/r^2}(t r^2)^{-d/2}e^{-r^2/(4tr^2)}
\mu(r^2 t) r^2\, dt\\
&\ge  c_8 r^{2-d}\int_0^1 t^{-d/2}e^{-1/(4t)} \mu(r^2 t)\, dt\ge
c_9
 r^{2-d}\int_0^1 t^{-d/2}e^{-1/(4t)} r^{-2}t^{-1}
\phi(r^{-2}t^{-1})\, dt\\
&\ge
c_{10} r^{-d}\int_0^1 t^{-d/2-1}e^{-1/(4t)}\phi(r^{-2})\, dt =
c_{11} r^{-d}\phi(r^{-2})\, ,
\end{align*}
where the last inequality follows because $r^{-2} t^{-1}\ge
r^{-2}$ and $\phi$ is increasing. \qed

Using \eqref{e:muatinfty} and \eqref{e:muat0}, we can easily show
(see \cite[Proposition 3.5]{KSV3} or \cite[Lemma 4.2]{RSV}) that

\begin{description}
\item{(1)} For any $M>0$, there exists $c=c(M, \phi)>0$ such that
\begin{equation}\label{H:1}
j(r)\le c j(2r), \qquad \forall r\in (0, M)\, .
\end{equation}
\item{(2)} There exists $c=c
(\phi)>0$ such that
\begin{equation}\label{H:2}
j(r)\le c j(r+1), \qquad \forall r>1.
\end{equation}
\end{description}

\section{Some results on symmetric L\'evy process in $\R$}\label{ss:1dsbm}
In this section we assume that $d=1$ and denote the process $X$ by $Z$. That is, $(Z_t,\P_x)$ is a
purely discontinuous symmetric L\'evy process in $\R$ such that
$$
\E_x\left[e^{i\xi\cdot(Z_t-Z_0)}\right]=e^{-t\Psi(|\theta|)},
\quad \quad \mbox{ for every } x\in \R \mbox{ and } \theta \in \R.
$$
We assume that \eqref{e:psi2} holds with a complete Bernstein function $\phi$ satisfying {\bf (H)}, that is, $\gamma^{-1}\phi(\theta^2)\le \Psi(|\theta|)\le \gamma\phi(\phi(\theta^2)$
for all $\theta\in \R$, but we do not assume the assumption \eqref{e:psi1} concerning the L\'evy measure of $Z$.
As a consequence of {\bf (H)}, \eqref{e:psi2} and \cite[Proposition 28.1]{Sa} we
know that for any $t>0$, $Z_t$ has a density $p_t(x, y)=p_t(y-x)$ which is smooth.

Let $\chi$ ($\kappa$, respectively)
 be the Laplace exponent of the
ladder height process of $Z$
($Y$, respectively).
It follows from \cite[Corollary 9.7]{Fris} that
\begin{equation}\label{e:formula4leoflh}
\chi(\lambda)
=\exp\left(\frac1\pi\int^{\infty}_0\frac{\log(\Psi(\lambda\theta))}
{1+\theta^2}d\theta \right), \quad
\kappa(\lambda)
=\exp\left(\frac1\pi\int^{\infty}_0\frac{\log(\phi(\lambda^2\theta^2))}
{1+\theta^2}d\theta \right)\, ,
\quad \forall \lambda>0.
\end{equation}
It follows immediately from these two equations and \eqref{e:psi2} that
$\gamma^{-1/2}\kappa(\lambda)\le \chi(\lambda)\le {\gamma}^{1/2}\kappa(\lambda)$,
i.e., that $\chi$ is comparable to $\kappa$.
From {\bf (H)} and
\cite[Propsoition 3.7]{KSV3} or \cite[Proposition 2.1]{KMR} we conclude that
the ladder height process of $Y$ has no drift and is not compound Poisson,
thus the ladder height process of $Z$ has no drift and is not compound Poisson.
Thus the process $Z$ does not creep upwards.
Since $Z$ is symmetric, we know that $Z$ also does not creep
downwards. Thus if, for any $a\in \R$, we define
$$
\tau_a=\inf\{t>0: Z_t<a\}, \quad \sigma_a=\inf\{t>0: Z_t\le a\},
$$
then we have
\begin{equation}\label{e:firstexittime}
\P_x(\tau_a=\sigma_a)=1, \quad x>a\, .
\end{equation}

Let $Z^{(0,\infty)}$ be the process $Z$ killed upon exiting $(0,\infty)$.
Since $Z$ has a smooth density, we can easily show that $Z^{(0,\infty)}$ has a density $p^{(0,\infty)}(t,x,y)$.
Let $G^{(0, \infty)}(x, y)=\int_0^{\infty}p^{(0,\infty)}(t,x,y)\, dt$ be the Green function of $Z^{(0,\infty)}$.
If we use $V$ to denote the potential measure of the ladder height process of
$Z$, then using the symmetry of $Z$ and \cite[Theorem 20, page 176]{Ber} we have that
for any $x\in (0, \infty)$ and any nonnegative function $f$ on $(0, \infty)$
\begin{equation}\label{e:gp4kphalfline}
\int^\infty_0f(y)G^{(0, \infty)}(x, y)dy=\int^\infty_0V(dy)\int^x_0V(dz)f(x+y-z).
\end{equation}
In the following, we will also use $V$ to denote the renewal function of the
ladder height process of $Z$: $V(t):=V((0,t))$.
For any $r>0$, let $G^{(0, r)}$ be the Green function
of $Z$ in  $(0, r)$.
Then we have the following result.

\begin{prop}\label{new1}
For all $r>0$ and all $x\in (0,r)$
$$
\int_0^r G^{(0,r)}(x,y)\, dy \le 2 V(r)\big(V(x) \wedge V(r-x)\big)\, .
$$
\end{prop}

\pf Since
$$
\int_0^r G^{(0,r)}(x, y)dy\le \int_0^r G^{(0,\infty)}(x, y)dy,
$$
we can apply \eqref{e:gp4kphalfline} with $f$ being the indicator function
of $(0, r)$ to immediately get the conclusion of the proposition.
\qed

The following result will play an important role in this paper.

\begin{prop}\label{new2}
There exists a constant $c=c(\gamma)>1$ such that for all $r>0$
$$
c^{-1} \frac1{\sqrt{\phi(r^{-2})}} \,\le\, V(r)\, \le\, c \frac1{\sqrt{\phi(r^{-2})}}.
$$
\end{prop}
\pf
The proof is a simple modification of
\cite[Theorem 4.4]{KMR}.
By \cite[Proposition 3.7]{KSV3} (or \cite[Proposition 2.1]{KMR}) we have that $c_1^{-1}\sqrt{\phi(\lambda)}\le \kappa(\lambda)\le c_1 \sqrt{\phi(\lambda)}$ for a constant $c_1>1$. Hence $c_2^{-1}\sqrt{\phi(\lambda)}\le \chi(\lambda)\le c_2 \sqrt{\phi(\lambda)}$, $c_2>1$, implying that
$$
c_3^{-1} \frac1{r \sqrt{\phi(r^{2})}} \le \sL V(r) \,\le\,
c_3 \frac1{r \sqrt{\phi(r^{2})}}
$$
(where $\sL V(r)$ denotes the Laplace transform of the \emph{function} $V$). The claim now follows by repeating the second part of the proof of \cite[Theorem 4.4]{KMR}.
\qed

\section{Poisson Kernel Estimates}

Recall that $Y$ is a subordinate Brownian motion in $\R^d$ with L\'evy exponent
$\phi(|\xi|^2)$,
$X$ is a purely discontinuous rotationally symmetric L\'evy process in $\R^d$
with L\'evy exponent
$\Psi(\xi)=\Psi(|\xi|)$
 and L\'evy density $J_X$, i.e.,
$$
\E_x\left[e^{i\xi\cdot(X_t-X_0)}\right]=e^{-t\Psi(|\xi|)},
\quad \quad \mbox{ for every } x\in \R^d \mbox{ and } \xi\in \R^d
$$
and
$\Psi(|\xi|)
= \int_{\R^d}(1-\cos(\xi\cdot y))J_X(y)dy.$
Recall that we assume that \eqref{e:psi1} holds.
As a consequence of {\bf (H)}, \eqref{e:psi2} and \cite[Proposition 28.1]{Sa} we
know that for any $t>0$, $X_t$ has a density $p_t(x, y)=p_t(y-x)$ which is smooth.

The infinitesimal generator $\L$ of $X$
is given by
\begin{equation}\label{3.1}
\L f(x)=\int_{\R^d}\left( f(x+y)-f(x)-y\cdot \nabla f(x)
{\bf 1}_{\{|y|\le1\}}
 \right)\, J_X(y)dy
\end{equation}
for $f\in C_b^2(\R^d)$. Moreover, for every $f\in C_b^2(\R^d)$,
$
f(X_t)-f(X_0)-\int_0^t \L f(X_s)\, ds
$
is a $\P_x$-martingale for every $x\in \R^d$.

First we record several inequalities  that will be needed in the remainder of the paper.
\begin{lemma}\label{l:l}
There exists a constant $c=c(\phi)>0$ such that
\begin{equation}\label{el6}
\int_0^{\lambda^{-1}} \phi (r^{-2})^{1/2}dr
 \,\le \, c \,
 \lambda^{-1}\phi(\lambda^{2})^{1/2},  \quad \forall \lambda \ge 1/R_0,
\end{equation}

 \begin{equation}\label{e:fsagd}
\lambda^2 \int_0^{\lambda^{-1}} r \phi (r^{-2})dr +
 \int_{\lambda^{-1}}^{R_0}r^{-1} \phi (r^{-2})dr \,\le\, c\,
\phi(\lambda^2), \quad \forall \lambda \ge 1/R_0,
\end{equation}
and
\begin{equation}\label{el3}
\lambda^2 \int_0^{\lambda^{-1}} r \phi (r^{-2})^{1/2}dr +
 \int_{\lambda^{-1}}^{R_0}r^{-1} \phi (r^{-2})^{1/2}dr \,\le\, c\,
\phi(\lambda^2)^{1/2}, \quad \forall \lambda \ge 1/R_0.
\end{equation}
\end{lemma}

\pf Assume $ \lambda \ge 1/R_0$. By {\rm (USC)}, $ \phi(r^{-2}) \le
c_1 r^{-2\delta_2}\lambda^{-2\delta_2}  \phi(\lambda^{2}) $ for $r \le \lambda^{-1}$.
On the other hand, by {\rm (LSC)}, $
 \phi(r^{-2}) \le  c_2 r^{-2\delta_1}\lambda^{-2\delta_1}  \phi(\lambda^{2}) $
for $\lambda^{-1}\le r\le R_0$.
Thus
\begin{align*}
\int_0^{\lambda^{-1}}  \phi (r^{-2})^{1/2}dr
\le  c_1^{1/2} \phi(\lambda^{2})^{1/2}\lambda^{-\delta_2}
\int_0^{\lambda^{-1}} r^{-\delta_2}   dr
&\le  c_3 \lambda^{-1}\phi(\lambda^{2})^{1/2} \frac{1}{1-\delta_2},
\end{align*}
\begin{align*}
&\lambda^2 \int_0^{\lambda^{-1}} r \phi (r^{-2})dr +
\int_{\lambda^{-1}}^{R_0} r^{-1} \phi (r^{-2})dr\\
&\le  c_4 \phi(\lambda^{2}) \left(\lambda^{2-2\delta_2}
\int_0^{\lambda^{-1}} r^{1-2\delta_2}   dr +\lambda^{-2\delta_1}
\int_{\lambda^{-1}}^{R_0} r^{-1-2\delta_1}  dr \right)\,
\le\,  c_5 \phi(\lambda^{2}) \left(\frac{1}{2(1-\delta_2)}  +
\frac{1}{2\delta_1} \right)
\end{align*}
and
\begin{align*}
&\lambda^2 \int_0^{\lambda^{-1}} r \phi (r^{-2})^{1/2}dr +
\int_{\lambda^{-1}}^{R_0} r^{-1} \phi (r^{-2})^{1/2}dr\\
&\le  c_6 \phi(\lambda^{2})^{1/2} \left(\lambda^{2-\delta_2}
\int_0^{\lambda^{-1}} r^{1-\delta_2}   dr +\lambda^{-\delta_1}
\int_{\lambda^{-1}}^{R_0} r^{-1-\delta_1}  dr \right)\,\le\,  c_7 \phi(\lambda^{2})^{1/2} \left(\frac{1}{2-\delta_2}  +
\frac{1}{\delta_1} \right).
\end{align*}
\qed

\begin{lemma}\label{l:lnew}
There exists a constant $c=c(\phi, \gamma)>0$ such that for every
$f\in C^2_b(\R^d)$ with $0\leq f \leq 1$,
$$
\L f_r(x) \le   c\, \phi(r^{-2}) \left( 2+\frac12 \sup_{y}\sum_{j,k} |(\partial^2/\partial
y_j\partial y_k) f(y)| \right) + b_0, \quad \text{for every } x \in \R^d, r  \le R_0
$$
where $f_r(y):=f(y/r)$ and $b_0:=2\int_{|z| > R_0} J_X(z)dz < \infty$.
\end{lemma}
\pf
Let $L_1=\sup_{y}\sum_{j,k} |(\partial^2/\partial
y_j\partial y_k) f(y)|$. Then $|f(z+y)-f(z) -y\cdot \nabla f(z)|\le
\frac{1}{2}L_1 |y|^2$. For $r\in (0,R_0]$, let $f_r(y)=f(y/r)$. Then the
following estimate is valid:
\begin{eqnarray*}
|f_r(z+y)-f_r(z) -y\cdot \nabla f_r(z){\bf 1}_{\{|y|\le r\}}| \le
\frac{L_1}{2}
\frac{|y|^2}{r^2}{\bf 1}_{\{|y|\le r\}} + 2 \cdot {\bf 1}_{\{|y|\ge r\}}\ .
\end{eqnarray*}
Now, by using
{\bf (H)},  \eqref{e:j-asymptotics} and \eqref{e:fsagd}, we get
\begin{eqnarray*}
|\L f_r(z)| &\le & \int_{\R^d}
|f_r(z+y)-f_r(z) -y\cdot \nabla f_r(z){\bf 1}_{\{|y|\le r\}}| \, J_X(y)dy \\
& \le & \frac{L_1}2\int_{\R^d} {\bf 1}_{\{|y|\le r\}}
\frac{|y|^2}{r^2}J_X(y)dy+2\int_{\R^d}{\bf 1}_{\{r \le |y|\le R_0\}} J_X(y)dy +2\int_{\R^d}{\bf 1}_{\{|y|\ge R_0\}} J_X(y)dy\\
& \le & \frac{\gamma L_1}2\int_{\R^d} {\bf 1}_{\{|y|\le r\}}
\frac{|y|^2}{r^2}j(|y|)dy+2\gamma\int_{\R^d}{\bf 1}_{\{r \le |y|\le R_0\}} j(|y|)dy +2\int_{\R^d}{\bf 1}_{\{|y|\ge R_0\}} J_X(y)dy\\
& \le &   c \phi(r^{-2}) \left( 2+\frac{L_1}2  \right) + 2\int_{\{|y|\ge R_0\}}J_X(y)dy \, ,
\end{eqnarray*}
where the constant $c$ is independent of $r\in (0,R_0]$.
\qed

For any open set $D$, we use $\tau_D$ to denote the first exit
time of $D$, i.e., $\tau_D=\inf\{t>0: \, X_t\notin D\}$.

Using Lemma \ref{l:lnew}, the proof of the next result is the same as those of
\cite[Lemmas 4.1 and 4.2]{KSV3}. Thus we skip the proof.

\begin{lemma}\label{L3.2}
There exists a constant $c=c(\phi,\gamma)>0$ such that for every $r\in (0,1]$, and every $x\in \R^d$,
$$
\inf_{z\in B(x,r/2)} \E_z \left[\tau_{B(x,r)} \right] \geq
\frac{c}{\phi((r/2)^{-2})}.
$$
\end{lemma}

The idea of the proof of the following proposition comes from \cite{Sz2}.
\begin{lemma}\label{l:tau}
There exists $c=c(\gamma)>0$ such that for any  $r\in (0, \infty)$ and $x_0 \in \R^d$,
\begin{eqnarray*}
\E_x[\tau_{B(x_0,r)}]\le  c\, (\phi(r^{-2})\phi((r-|x-x_0|)^{-2}))^{-1/2}
\qquad x\in B(x_0, r).
\end{eqnarray*}
\end{lemma}

\pf Without loss of generality, we may assume that $x_0=0$. We fix
$x\neq 0$ and  put $Z_t=\frac{X_t\cdot x}{|x|}$. Then, using the fact that $\Psi$ is a radial function,
$Z_t$ is a L\'evy process on $\R$ with
$$
\E[e^{i\theta Z_t}]=\E(e^{i\theta\frac{x}{|x|}\cdot X_t})
=e^{-t \Psi(\theta\frac{x}{|x|})}
=e^{-t \Psi(\theta)},
\qquad \theta\in \R.
$$
Clearly, $\gamma^{-1}\phi(\theta^2)\le \Psi(\theta)\le \gamma\phi(\theta^2)$.
Thus $Z_t$ is of the type of one-dimensional symmetric L\'evy
processes studied in Section \ref{ss:1dsbm}.

It is easy to see that, if
$X_t\in B(0, r)$, then $|Z_t|<r$, hence
$
\E_x[\tau_{B(0, r)}]\le \E_{|x|}[\tilde \tau],
$
where $\tilde \tau=\inf\{t>0: |Z_t|\ge r\}$.
Thus, applying  Proposition \ref{new1}, we
obtain
$
\E_x[\tau_{B(0, r)}] \le 2 V(2r) V(r-|x|).
$
Now, by Proposition \ref{new2} and {\bf (H)}, we
have proved the lemma.
\qed

 We now recall the definition
of harmonic functions with respect to $X$.

\begin{defn}\label{def:har1}
Let $D$ be an open subset of $\R^d$.
A function $u$ defined on $\R^d$ is said to be

\begin{description}
\item{\rm{(1)}}  harmonic in $D$ with respect to $X$ if
$$
\E_x\left[|u(X_{\tau_{B}})|\right] <\infty
\quad \hbox{ and } \quad
u(x)= \E_x\left[u(X_{\tau_{B}})\right],
\qquad x\in B,
$$
for every open set $B$ whose closure is a compact
subset of $D$;

\item{\rm{(2)}}
regular harmonic in $D$ with respect to $X$ if it is harmonic in $D$
with respect to $X$ and
for each $x \in D$,
$$
u(x)= \E_x\left[u(X_{\tau_{D}})\right].
$$
\end{description}
\end{defn}

Since our $X$ satisfies \cite[(1.6), {\bf (UJS)}]{CKK}, by \cite[Theorem 1.4]{CKK} and  using the standard chain argument one have  the following
form of Harnack inequality.

\begin{thm}\label{hi}
For every $a \in (0,1)$, there exists $c=c(a, \phi, \gamma)>0$ such that
for every $r \in (0, 1)$, $x_0 \in {\mathbb R}^d$, and any function $u$ which is nonnegative on ${\mathbb R}^d$
and harmonic with respect to $X$ in $B(x_0, r)$, we have
$$
u(x)\,\le \,c\, u(y), \quad \textrm{for all }x, y\in B(x_0, ar)\, .
$$
\end{thm}

Given  an open set $D\subset \R^d$, we define
$X^D_t(\omega)=X_t(\omega)$ if $t< \tau_D(\omega)$ and
$X^D_t(\omega)=\partial$ if $t\geq  \tau_D(\omega)$, where
$\partial$ is a cemetery state.
A subset $D$ of $\R^d$ is said to be Greenian (for $X$) if $X^{D}$ is transient. When
$d\ge 3$, any non-empty open set $D\subset\bR^d$ is
Greenian. An open set $D\subset \bR^d$ is Greenian if
and only if $D^c$
 is non-polar for $X$ (or equivalently, has positive capacity with respect to $X$).
 In particular, every bounded open
set is Greenian.

Since $X$ has a smooth density, using the strong Markov property,  it is standard to show that for every Greeninan open set $D$, $X^D_t$ has a density $p_D(t,x,y)$. For any Greeninan open  set $D$ in $\R^d$ let $G_D(x,y)=\int_0^{\infty}p_D(t,x,y)$ be the Green function of $X^D$.
Using the L\'{e}vy system for $X$, we know that for every
Greeninan open subset $D$ and  every $f \ge 0$ and $x \in D$,

\begin{equation}\label{newls}
\E_x\left[f(X_{\tau_D});\,X_{\tau_D-} \not= X_{\tau_D}  \right]
=  \int_{\overline{D}^c} \int_{D}
G_D(x,z) J_X(z-y) dz f(y)dy.
\end{equation}
We define the Poisson kernel
\begin{equation}\label{PK}
K_D(x,y)\,:=  \int_{D}
G_D(x,z) J_X(z-y) dz, \qquad (x,y) \in
D \times  {\overline{D}}^c.
\end{equation}
Thus \eqref{newls} can be simply written as
$$
\E_x\left[f(X_{\tau_D});\,X_{\tau_D-} \not= X_{\tau_D}  \right]
=\int_{\overline{D}^c} K_D(x,y)f(y)dy.
$$
Using  continuity of $J_X$, one can easily check
that $K_D(x,\cdot)$ is continuous on $\overline{D}^c$ for every $x\in D$.

\begin{prop}\label{p:Poisson1}
There exists $c_1=c_1
(\phi,\gamma)>0$ and $c_2=c_2
(\phi, \gamma)>0$ such that for every $r \in (0, 1]$ and $x_0 \in \R^d$,
\begin{eqnarray}
K_{B(x_0,r)}(x,y) \,&\le &\, c_1 \, j(|y-x_0|-r) \left(\phi(r^{-2})\phi((r-|x-x_0|)^{-2})\right)^{-1/2} \label{P1}\\
 &\le &\, c_1 \, j(|y-x_0|-r) \phi(r^{-2})^{-1}\ \label{P1-worse}
\end{eqnarray}
for all $(x,y) \in B(x_0,r)\times \overline{B(x_0,r)}^c$ and
\begin{equation}\label{P2}
K_{B(x_0, r)}(x_0, y) \,\ge\, c_2\, j(|y-x_0|)
\phi(r^{-2})^{-1}, \qquad \textrm{ for all } y \in \overline{B(x_0, r)}^c.
\end{equation}
\end{prop}

\pf
Using \eqref{e:psi1} and \eqref{H:1}--\eqref{H:2}, the proof of \eqref{P1} and \eqref{P2} is exactly the same
as that of \cite[Proposition 4.10]{KSV3}
(using {\bf (H)}), while \eqref{P1-worse} follows from \eqref{P1} and the fact that $\phi$ is increasing.
\qed

Using Theorem \ref{hi} and the continuity of
$K_{B(x_0,r)}(x,\cdot)$  on $\overline{B(x_0,r)}^{\, c}$ for every $x\in D$,
 the proof of the next result is the same as that of
\cite[Proposition 1.4.11]{KSV3}. So we omit the proof.
\begin{prop}\label{p:Poisson2}
For every $a \in (0,1)$, there exists $c=c(\phi, \gamma, a)>0$ such that for every
$r \in (0,1]$, $x_0 \in {\mathbb R}^d$ and $x_1, x_2 \in B(x_0, ar)$,
$$
K_{B(x_0,r)}(x_1,y) \,\le\, c  K_{B(x_0,r)}(x_2,y), \qquad y \in
\overline{B(x_0,r)}^{\, c}\, .
$$
\end{prop}

\begin{prop}\label{p:Poisson3}
For every $a \in (0,1)$, there exists $c=c(\phi, \gamma, a)>0$ such that for every
$r \in (0, 1]$ and $x_0 \in \R^d$,
$$
K_{B(x_0,r)}(x,y)\le c\ r^{-d} \left(\frac{\phi((|y-x_0|-r)^{-2})}{\phi(r^{-2})}\right)^{1/2}
$$
for all $x\in B(x_0,ar)$ and all $y$ such that $r<|x_0-y|<2r$.
\end{prop}
\pf
By Proposition \ref{p:Poisson2},
$$
K_{B(x_0,r)}(x,y) \le \frac{c_1}{r^d} \int_{B(x_0, a r)}
K_{B(x_0,r)}(w,y) dw
$$
for some constant $c_1=c_1(\phi, \gamma, a)>0$. Thus from Lemma \ref{l:tau}, \eqref{PK} and
Theorem \ref{t:Jorigin} we have that
\begin{eqnarray*}
K_{B(x_0,r)}(x,y)&\le& \frac{
c_1}{r^d}\int_{B(x_0, r)}\int_{B(x_0,r)} G_{B(x_0,r)}(w,z)J_X(z-y) dz dw \\
&=& \frac{c_1}{r^d}\int_{B(x_0, r)} \E_z[\tau_{B(x_0,r)}]J_X(z-y) dz\\
& \le& \frac{c_2}{r^d(\phi(r^{-2}))^{1/2}}
\int_{B(x_0, r)}\frac{\phi(|z-y|^{-2})}{(\phi((r-|z-x_0|)^{-2}))^{1/2}} |z-y|^{-d} dz
\end{eqnarray*}
for some constant $c_2=c_2(\phi, \gamma, a)>0$.
Since $r-|z-x_0| \le |y-z|$, we have
\begin{eqnarray*}
K_{B(x_0,r)}(x,y) & \le&
\frac{c_2}{r^d(\phi(r^{-2}))^{1/2}}
\int_{B(x_0, r)}\frac{(\phi(|z-y|^{-2}))^{1/2}}{|z-y|^{d}} dz\\
& \le&   \frac{c_2}{r^d(\phi(r^{-2}))^{1/2}}\int_{B(y, 3r) \setminus B(y,
|y-x_0|-r)}
\frac{(\phi(|z-y|^{-2}))^{1/2}}{|z-y|^{d}} dz\\
&\le & \frac{c_3}{r^d(\phi(r^{-2}))^{1/2}}
\int_{|y-x_0|-r}^{3r } \frac{\phi(s^{-2})^{1/2}}
{s}ds.
\end{eqnarray*}
Now, using \eqref{el3} in the
last integral (considering the cases $r < R_0/3$ and $1 \ge r \ge R_0/3$ separately),
we arrive at the conclusion of the proposition.
 \qed

\begin{lemma}\label{l2.1}
For every $a \in (0, 1)$, there exists a positive constant $c=c(\phi, \gamma, a)>0$ such that
for any $r\in (0, 1)$ and any open set $D$ with $D\subset B(0, r)$ we have
$$
{\P}_x\left(X_{\tau_D} \in B(0, r)^c\right) \,\le\, c\,
\phi(r^{-2})\int_D G_D(x,y)dy, \qquad x \in D\cap B(0, ar)\, .
$$
\end{lemma}
\pf
The proof of the lemma is similar to that of \cite[Lemma 4.15]{KSV3}.
Transience was used in the proof \cite[Lemma 4.15]{KSV3} in order to derive equation  \cite[(4.29)]{KSV3}.
By noting that \cite[(4.29)]{KSV3} in the proof of  \cite[Lemma 4.15]{KSV3} follows immediately from Dynkin's formula,
using our Lemma \ref{l:lnew}, we can follow the rest of the
proof \cite[Lemma 4.15]{KSV3} (which does not use transience) to get the conclusion of the lemma here. We omit the details.
\qed

\section{Uniform Boundary Harnack Principle}

In this section, we give the proof of the main result of this paper.
Let  $A(x, a,b):=\{ y \in \R^d: a \le |y-x| <b  \}.$

\begin{lemma}\label{l:k_B}
For every $p \in (0,1)$, there exists
$c=c(\phi, \gamma, p)>0$ such that for every $r \in (0, 1)$,
$$
\int_{r(1+p)/2}^{|y|}
  K_{B(0,s)}(x,y)ds \,\le\,
c\,\frac{r}{\phi(r^{-2})}  j(|y|) \qquad \forall
x\in  B(0, pr),\,  y\in A(0, r(1+p)/2, r).$$
\end{lemma}
\pf
Let $0<p<1$ and $q=(1+p)/2$.
Note that the functions $r\mapsto r^{-d+1}$ and $r\mapsto r^{-1}(\phi(r^{-2}))^{-1/2}$ are decreasing, see \eqref{e:uv}.
Using Proposition \ref{p:Poisson3} we get
\begin{eqnarray*}
\int_{qr}^{|y|}
  K_{B(0,s)}(x,y)ds &\le& c_1\,    \int_{qr}^{|y|}\frac{s^{-d}}{(\phi(s^{-2}))^{1/2}}
(\phi((|y|-s)^{-2}))^{1/2}ds  \\
&\le & c_2\,   \frac{r^{-d}}{(\phi((qr)^{-2}))^{1/2}} \int_{qr}^{|y|}
(\phi((|y|-s)^{-2}))^{1/2}ds
\end{eqnarray*}
for some constants $c_1 (p, \phi)>0$ and $c_2 (p, \phi)>0$.
Note that  by \eqref{el6} (considering the cases $|y|-qr< R_0$ and $1 \ge |y|-qr \ge R_0$ separately)
and the fact that $r \mapsto r (\phi(r^{-2}))^{1/2}$ is increasing
$$
  \int_{qr}^{|y|}
(\phi((|y|-s)^{-2}))^{1/2}
ds = \int_{0}^{|y|-qr}
(\phi(s^{-2}))^{1/2}
ds \le c_3 (|y|-qr) (\phi((|y|-qr)^{-2}))^{1/2} \le c_3 r (\phi(r^{-2}))^{1/2}$$
for some constant $c_3 >0$. Thus, by {\bf (H)}, Theorem \ref{t:Jorigin} and the fact that $r \to j (r)$ is decreasing,
we have
\begin{eqnarray*}
\int_{qr}^{|y|}
  K_{B(0,s)}(x,y)ds \le \frac{c_4}{r^{d-1}} \le c_5 \frac{r}{ \phi(r^{-2})} j(r) \le c_6 \frac{r}{ \phi(r^{-2})} j(|y|).
\end{eqnarray*}
\qed

From the strong Markov property, it is well known and easy to see that
for every Greenian open sets $U$ and $D$ with $U \subset D$,
$
G_D(x,y)= G_U(x,y) + \E_x\left[   G_D(X_{\tau_U}, y)\right]
$ for every $(x,y) \in \R^d \times \R^d.
$
Thus,
for every Greenian open sets $U$ and $D$ with $U \subset D$,
\beq\label{e:KDU}
K_D(x,z)= K_U(x,z) + \E_x\left[   K_D(X_{\tau_U}, z)\right], \qquad (x,z) \in
U \times \overline{D}^c
\eeq
and
\beq\label{e:EDU}
\E_x[\tau_D]= \E_x[\tau_U] + \E_x\left[   \E_{X_{\tau_U}}[\tau_D]\right], \qquad x \in U.
\eeq

\begin{lemma}\label{lK1_1}
For every $p \in (0,1)$, there exists $c=
c(\phi, \gamma, p)>0$ such that for every $r \in (0, 1)$,
for every $z_0 \in \R^d$, $U \subset B(z_0,r)$ and for any  $ (x,y) \in (U\cap B(z_0, pr)) \times B(z_0, r)^c$,
\begin{eqnarray*}
K_U(x, y)
&\le&
c\,\frac{1}{\phi(r^{-2})} \left(
\int_{U\setminus B(z_0, (1+p)r/2)}  j(|z-z_0|)  K_U(z,y)dz +j(|y-z_0|)
\right).
\end{eqnarray*}
\end{lemma}
\pf
Without loss of generality, we assume $z_0=0$.
Let $0 <p<1$, $q_1:=(1+p)/2$ and $q_2:=(3+2p)/5. $
For every $s \in [q_1 r,q_2 r]$ and $x\in U \cap B(0,pr)$,
by \eqref{e:KDU}  we have
\begin{eqnarray*}
 K_U(x, y)
   &=&
   \E_x[K_U(X_{\tau_{U \cap B(0,s)}},y)]  + K_{U \cap B(0,s)}(x, y)\\
     &=&\int_{U\setminus B(0,s)}K_U(z,y)
  K_{U \cap B(0,s)}(x,z)dz
+
K_{U \cap B(0,s)}(x,y)\\
  &\le&\int_{U\setminus B(0,s)}K_U(z,y)
  K_{B(0,s)}(x,z)dz
+
K_{B(0,s)}(x,y).
\end{eqnarray*}
Thus
\begin{align*}
 K_U(x, y)&\le \frac{1}{r(q_2 - q_1)}\int_{q_1r}^{q_2r} \int_{U\setminus B(0,s)} K_{B(0,s)}(x,z) K_U(z,y)  dzds
 + \frac{1}{r(q_2 - q_1)}\int_{q_1r}^{q_2r}K_{B(0,s)}(x,y)ds\\
 &=:I+II.
\end{align*}
By Tonelli's theorem, we have
\begin{eqnarray*}
I  &=&\frac{10}{r(1 - p)} \int_{q_1r}^{q_2r} \int_{\{z \in U; |z|\ge q_1r\}}1_{\{|z|\ge s\}} K_{B(0,s)}(x,z) K_U(z,y)  dzds\\
&\le & \frac{10}{r(1 - p)}
\int_{(U\setminus B(0, q_1r)) }\left(\int_{q_1r}^{|z|}
  K_{B(0,s)}(x,z)ds \right) K_U(z,y) dz.
\end{eqnarray*}
Applying Lemma \ref{l:k_B} to the inner integral above, we get that
\bee\label{e:w1}
I\leq
\frac{c_1}{(1 - p)}\frac{1}{\phi(r^{-2})}
\int_{(U\setminus B(0, q_1r))}  j( |z|)  K_U(z,y) dz.
\eee

One the other hand,  for any $s\in [q_1r, q_2r]$, by Proposition \ref{p:Poisson1},
$$
K_{B(0,s)}(x,y)
\le c_2 j(|y|-s) \frac{1}
{(\phi(s^{-2}))^{1/2}}\frac{1}
{(\phi((s-|x|)^{-2}))^{1/2}}.
$$
When $y \in A(0,  r , 4)$ we have $(1-q_2)|y|\le |y|-s $,
while when $|y|\ge 4$ we have $|y|-s\ge |y|-1$. Since $
s-|x|\le s \le q_2{r}$, we have by the monotonicity of $j$,
$$
j(|y|- s )
\frac{ 1}{(\phi( s ^{-2}))^{1/2}}
\frac{1}
{(\phi(( s -|x|)^{-2}))^{1/2}}
  \,\le\, c_3 j((1-q_2)|y|)
\frac{1}{\phi(r^{-2})}
, \quad y \in A(0,  r , 4)
$$
and
$$
j(|y|- s ) \frac{1}{(\phi( s ^{-2}))^{1/2}}
\frac{1} {(\phi(( s -|x|)^{-2}))^{1/2}}
 \,\le\, c_3 j(|y|-1)
\frac{1}{\phi(r^{-2})} , \quad |y|\ge 4
$$
for some constant $c_3>0$. Thus by applying \eqref{H:1} and
\eqref{H:2}, we get
\bee\label{e:w2}
II\le c_4(1 - p)^{-1}\int_{q_1r}^{q_2r} j(|y|- s )
\frac{ 1}{(\phi( s ^{-2}))^{1/2}}
\frac{1}
{(\phi(( s -|x|)^{-2}))^{1/2}} ds \le c_5 j(|y|)
\frac{1}{\phi(r^{-2})}.
\eee
Combining \eqref{e:w1}-\eqref{e:w2}, we conclude that
$$
K_U(x,y)\leq
c_6 \frac{1}{\phi(r^{-2})}
\int_{(U\setminus B(0, q_1r))}  j( |z|)  K_U(z,y) dz
+c_6 j(|y|)
\frac{1}{\phi(r^{-2})}.
$$
\qed

Note that, since $X$ satisfies the hypothesis ${\bf H}$ in \cite{Sz1},
by \cite[Theorem 1]{Sz1}, if $V$ is a Lipschitz open set and $U \subset V$
\beq\label{e:Lip}
\P_x(X_{\tau_U} \in \partial V)=0 \qquad \text{ and } \qquad \P_x(X_{\tau_U} \in dz) = K_U(x,z)dz \quad\text{on }V^c.
\eeq

\begin{lemma}\label{lK1}
For every $p \in (0,1)$, there exists $c=
c(\phi, \gamma, p)>0$ such that for every
$r \in (0, 1)$, for every $z_0 \in \R^d$, $U \subset B(z_0,r)$ and  any nonnegative function $u$ in $\RR^d$
which is regular harmonic in $U$ with respect to $X$ and vanishes in
$U^c \cap B(z_0, r)$ we have
\begin{eqnarray*}
u(x)
&\le&
c\,\frac{1}{\phi(r^{-2})}
\int_{(U\setminus B(z_0, (1+p)r/2)) \cup B(z_0,r)^c}  j(
|y-z_0|) u(y)dy,\quad
x \in U \cap B(z_0, pr).
\end{eqnarray*}
\end{lemma}
\pf
Without loss of generality, we assume $z_0=0$.
Let $0 <p<1$ and set $
q=(1+p)/2$.
Note that the part of boundary of $U$ belonging to $U^c \cap B(z_0, r)$ needs not be Lipschitz,
but here $u$ vanishes. The other part of boundary of $U$ is a part of the boundary of the ball $B(z_0,r)$.
Thus, since $u$ is regular harmonic in $U$ with respect to $X$ and
vanishes in $U^c \cap B(z_0, r)$,
by  Lemma \ref{lK1_1} and \eqref{e:Lip} we have
\begin{eqnarray*}
u(x)&= &\E_x[u(X_{\tau_U})]=\int_{U^c} K_U(x,y)u(y)dy \\
  &\le& c \frac{1}{\phi(r^{-2})}\left( \int_{U^c}\int_{U\setminus B(0,
qr)} j( |z|)  K_U(z,y) dz u(y)dy  +
  \int_{B(0,r)^c} j(|y|)  u(y)dy\right).
\end{eqnarray*}
Since $\int_{U^c} K_U(z,y)u(y)dy=u(z)$ on $U\setminus B(0,
qr)$, by Tonelli's theorem, we have
\begin{eqnarray*}
u(x)\,\le\,c \frac{1}{\phi(r^{-2})}\left( \int_{U\setminus B(0,
qr)} j( |z|)   u(z) dz +
  \int_{B(0,r)^c} j(|y|)   u(y)dy\right).
\end{eqnarray*}
\qed

\begin{lemma}\label{lK2_1}
There exists $C=C(\phi, \gamma
)>1$ such that for every $r \in (0, 1)$,
for every $z_0 \in \R^d$, $U \subset B(z_0,r)$ and for any  $ (x,y) \in U \cap B(z_0, r/2) \times
(B(z_0, r)^c \cap {\overline{U}}^c)$,
\begin{eqnarray*}
\lefteqn{C^{-1}\, \E_x[\tau_{U}]
 \left(\int_{U\setminus B(z_0, r/2)}  j(|z-z_0|)  K_U(z,y)dz  +j(|y-z_0|) \right)} \\
&\le &K_U(x,y) \le
C\,  \E_x[\tau_{U}] \left(\int_{U\setminus B(z_0, r/2)}  j(|z-z_0|)  K_U(z,y)dz  +j(|y-z_0|) \right).
 \end{eqnarray*}
 \end{lemma}
\pf
Without loss of generality, we assume $z_0=0$. Fix $r \in (0, 1)$ and let $B:=B(0,r)$,
$U_1:=U \cap B(0, \frac12 r)$, $U_2:=U \cap B(0, \frac23 r)$
 and $U_3:=U \cap B(0, \frac34 r)$.
Let $x\in U\cap B(0, r/2)$,
$y\in B(0,r)^c    \cap {\overline{U}}^c$.
By \eqref{e:KDU},
\begin{eqnarray*}
K_{U}(x,y)
&=&\E_x[K_U(X_{\tau_{U_2}},y)]  + K_{U_2}(x, y)\\
  &=&\int_{U_3\setminus U_2}K_U(z,y)\P_x(X_{\tau_{U_2}} \in dz) +\int_{U\setminus U_3}K_U(z,y)K_{U_2}(x,z)dz+  K_{U_2}(x, y) \\
  &=&\int_{U_3\setminus U_2}K_U(z,y)\P_x(X_{\tau_{U_2}} \in dz)+\int_{U\setminus U_3}K_U(z,y)\int_{U_2}G_{U_2}(x,w)j(|z-w|)dw dz\\
  &&\quad+
  \int_{U_2}G_{U_2}(x,w)j(|y-w|)dw=:I+II+III.\\
 \end{eqnarray*}
From Lemmas \ref{l2.1} and \ref{lK1_1}, we see that there exist
$c_i=c_i(\phi, \gamma)>0, i=1, 2,$ such that $I$ is less than or equal to
  \beq \label{e:k1}
 c_1\left(\sup_{z \in U_3}K_U(z,y)\right)\phi(r^{-2})\E_x[\tau_{U_2}]   \,\le\, c_2\E_x[\tau_{U_2}]
 \left(\int_{U\setminus U_3}   j( |z|)  K_U(z,y) dz  + j(|y|) \right). \eeq

 Note that if $4>|z| \ge \frac34 r$ and $|w| < \frac23 r$, then
 $$
 \frac19|z| = |z|-\frac89 |y| \le |z|-\frac23r \le |z|-|w| \le |z-w| \le |z|+|w| \le |z|+\frac23 r \le2|z|.
 $$
 Thus, since $j$ is monotonely decreasing, by \eqref{H:1}
 $$
 c_3^{-1} j(|z|) \,\le \, j(2 |z|) \, \le\, j(|z-w|)\, \le\,  j(\frac19 |z|)\,\le\, c_3 j(|z|),
 \qquad \text{ if }  4>|z| \ge \frac34 r, |w| < \frac23 r
 $$ for some constant $c_3>0$.
 If $4 \le |z|$ and $|w| < \frac23 r$, then $|z|-\frac23r \le |z-w| \le |z|+\frac23 r$ and by
 \eqref{H:2}
$$
c_4^{-1} j(|z|) \,\le \, c_4^{-1} j(|z|+\frac23 r-1) \,\le \, j(|z|+\frac23 r) \, \le\, j(|z-w|)\, \le\,  j(|z|-\frac23 r)\,\le\, c_4 j(|z|-\frac23 r+1)\,\le\, c_4 j(|z|),
$$
for some constant $c_4>0$. Thus
there exists $c_5=c_5(\phi, \gamma)>1$ such that
  \beq \label{e:k2}
c_5^{-1}  \E_x[\tau_{U_2}]
  \int_{U\setminus U_3}  j(|z|)  K_{U}(z,y)dz \le II \le c_5\E_x[\tau_{U_2}]\int_{U\setminus U_3 } j(|z|)  K_{U}(z,y)dz
\eeq
and
 \beq \label{e:k2_1}
 c_5^{-1} \E_x[\tau_{U_2}]j(|y|) \le III  \le c_5 \E_x[\tau_{U_2}]j(|y|)\, .
\eeq
Now the upper bound follows from \eqref{e:k1}--\eqref{e:k2_1}. To prove the lower bound we can neglect $I$. Further,
by using Lemmas \ref{l:tau} and \ref{l2.1} in the third line, from  \eqref{e:EDU} we get
  \begin{eqnarray*}
\E_x[\tau_{U}]&= &\E_x[\tau_{U_2}] + \E_x\left[   \E_{X_{\tau_{U_2}}}[\tau_{U}]\right]\\
&\le &\E_x[\tau_{U_2}] + \left(\sup_{z \in U}\E_z[\tau_{U}]\right)  \P_x\big(X_{\tau_{U_2}} \in B(0, 2r/3 )^c\big) \\
&\le &\E_x[\tau_{U_2}] + c_6\phi(r^{-2})^{-1}\,  c_7\phi((2r/3)^{-2})\E_x[\tau_{U_2}]\,\le \, c_8 \E_x[\tau_{U_2}]
 \end{eqnarray*}
 for some constants $c_8 >0$.
In the last inequality above we have used {\bf (H)}.
 Since
\begin{eqnarray*}
  \int_{U\setminus U_1}  j(|z|)  K_U(z,y)dz
  &=&\int_{U\setminus U_3}  j(|z|)  K_U(z,y)dz+\int_{U_3 \setminus U_1}
  j(|z|)  K_U(z,y)dz \\
    &\le &\int_{U\setminus U_3}  j(|z|)  K_U(z,y)dz +\left(\sup_{z \in U_3} K_U(z,y)
    \right) \int_{A(0, r/2, 3r/4)}
   j(|y|)  dy,
    \end{eqnarray*}
  by Theorem \ref{t:Jorigin} and Lemma
  \ref{lK1_1},
\begin{eqnarray*}
   \int_{U\setminus U_1}  j(|z|)  K_U(z,y)dz
  & \le& \left(1+ \frac{c_9}{\phi(r^{-2})}\int_{r/2}^{3r/4}  s^{-1} \phi (s^{-2})  ds \right)
 \left(\int_{U\setminus U_3}  j(|z|)  K_U(z,y)dz + j(|y|)\right).
\end{eqnarray*}
  Applying \eqref{e:fsagd} (considering the cases $r < 4R_0/3$ and $1 \ge r \ge 4R_0/3$ separately) and {\bf (H)}, we obtain
  \begin{eqnarray}
 \int_{U\setminus U_1}  j(|z|)  K_U(z,y)dz
  \le  c_{10} \left(\int_{U\setminus U_3}  j(|z|)  K_U(z,y)dz + j(|y|)\right). \label{e:k4}
  \end{eqnarray}

  Combining \eqref{e:k2}-\eqref{e:k4}, we have proved the lower bound.
\qed

\begin{lemma}\label{lK2}

For every $z_0 \in \R^d$, every open set $U \subset B(z_0,r)$ and for any nonnegative function $u$ in $\RR^d$ which is regular harmonic in $U$ with respect to $X$ and vanishes a.e.~in $U^c \cap B(z_0, r)$ it holds that
\begin{eqnarray*}
C^{-1}\E_x[\tau_{U}] \int_{B(z_0, r/2)^c}  j(|y-z_0|)  u(y)dy
\le u(x)\le
C \E_x[\tau_{U}] \int_{B(z_0, r/2)^c} j(|y-z_0|)  u(y)dy
\end{eqnarray*}
for every  $ x \in U \cap B(z_0, r/2)$ (where $
C$ is the constant from Lemma \ref{lK2_1}).
\end{lemma}
\pf
Without loss of generality we may take $z_0=0$. By the argument in the proof of Lemma \ref{lK1} and by the assumption that $u$ vanishes a.e.~on $U^c \cap B(0, r)$ we have that
$$
u(x)=\int_{U^c}K_U(x,y)u(y)\, dy =\int_{B(z_0,r)^c} K_U(x,y)u(y)\, dy\, .
$$
Now the claim follows from Lemma \ref{lK2_1}.
Indeed, by Tonelli's theorem we get
\begin{eqnarray*}
u(x)&\le &
C \, \E_x[\tau_U]\left(\int_{U\setminus B(0,r/2)}j(|z|)
\left(\int_{B(0,r)^c}K_U(z,y)u(y)\, dy\right)dz+\int_{B(0,r)^c}j(|y|)u(y)\, dy\right)\\
&=&
C \, \E_x[\tau_U] \left(\int \int_{U\setminus B(0,r/2)}
j(|u|)u(z)\, dz +\int_{B(0,r)^c}j(|z|)u(z)\, dz\right)\\
&=&
C \, \E_x[\tau_U]\int_{B(0,r/2)}j(|u|)u(z)\, dz,
\end{eqnarray*}
where for the last line we used that $u$ vanishes a.e.~on $U^c\cap B(0,r)$.
The lower bound follows in the same way. \qed

We remark that in the statements of Lemmas \ref{l:k_B}--\ref{lK2}, by using the assumption \eqref{e:psi1}, we could have replaced the density $j$ with the density $J_X$ of the process $X$ (with a different constant). We will do this in the next corollary which gives an approximate factorization of the Poisson kernel. It is an immediate consequence of the last two lemmas.

\begin{corollary}\label{c:approx-factor}
Let $z_0\in \R^d$, $D\subset \R^d$ be
Greenian open set and denote $U:=D\cap B(z_0,r)$.
Then for every $r\in (0,1)$ and all $(x,y)\in (D\cap  B(z_0,r/2))\times (D^c\cap B(z_0,r)^c)$  it holds that
\begin{equation}\label{e:approx-factor}
C^{-1} \E_x[\tau_U] A(y) \le  K_D(x,y) \le
C\, \E_x[\tau_U] A(y)\, ,
\end{equation}
where
$$
A(y):=\int_{U\setminus B(z_0,r/2)}\big(J_X(z-z_0)K_U(z,y)\, dz+J_X(y-z_0)\big)+\int_{B(z_0,r/2)^c} J_X(z-z_0)\E_z\left[K_D(X_{\tau_U},y)\right]\, dz \, .
$$
\end{corollary}
\pf Without loss of generality, we assume $z_0=0$ and $D \cap B(0, r/2) \not= \emptyset$. We first note that by \eqref{e:KDU} and \eqref{e:Lip}, for every $(x,y) \in (D \cap B(0, r)) \times (D^c \cap B(0, r)^c)$,
$$
K_D(x, y) =  K_U(x,y) + \E_x\left[   K_D(X_{\tau_U}, y)\right].
$$
The function $x\mapsto \E_x\left[   K_D(X_{\tau_U}, y)\right]$ is regular harmonic in $U$
with respect to $X$ and vanishes a.e.~in $U^c\cap B(0,r)$. By using Lemma \ref{lK2}
for this function, and Lemma \ref{lK2_1} for $K_U(x,y)$ we immediately obtain required inequalities. \qed

We can now easily prove Theorem \ref{UBHP}.

\noindent
{\bf Proof of Theorem \ref{UBHP}.} (i) This follows immediately from Lemma \ref{lK2} with
$c:=C^4$.

\noindent (ii) Let $x_1,x_2\in D\cap B(z_0,r/2)$, $y_1,y_2\in D^c\cap B(z_0,r)^c$ and let $U:=D\cap B(z_0,r)$. Then by \eqref{e:approx-factor}
\begin{eqnarray*}
K_D(x_1,y_1) K_D(x_2, y_2)&\le & \left(
C \E_{x_1}[\tau_U]A(y_1)\right)\left(
C \E_{x_2}[\tau_U]A(y_2)\right)\\
&=&\left(
C \E_{x_1}[\tau_U]A(y_2)\right)\left(
C \E_{x_2}[\tau_U]A(y_1)\right)\\
&\le &
C^4 K_D(x_1,y_2) K_D(x_1, y_2)\, .
\end{eqnarray*}
The lower bound is proved in the same way. \qed

\section{Remarks on {\bf (H)}}

In this section we point out the relationship between the assumption
{\bf (H)} and the class $OR$ of $O$-regularly varying functions,
and sketch the construction of a complete Bernstein function $\phi$
which satisfies {\bf (H)} but not the assumption in \cite{KSV3} that $\phi$
is comparable to a regularly varying function. Using the idea in the construction below, one
can come up with a complete Bernstein function that is bounded between any two regularly varying
complete Bernstein functions.

It follows from the definitions on \cite[page 65 and page 68]{BGT} and
\cite[Proposition 2.2.1]{BGT} that the assumption {\bf (H)} is equivalent to
that $\phi$ is in $OR$ with its Matuszewska indices contained in $(0, 1)$.

Here is a sketch of the construction.
For $x\in (0, 2]$, define
$$
f(x)=x^{1/2}.
$$
Then we define
$$
f(x)=x^{1/3}+f(2)-2^{1/3}, \qquad x\in (2, a_1]
$$
for some large constant $a_1>2$. The constant $a_1$ is chosen so that for large
values of $x$ in $(2, a_1]$, the function $f$ behaves like $x^{1/3}$, that is
$f(\lambda x)/f(x)$ is close to $\lambda^{1/3}$ uniformly for $\lambda\in [1, 2]$.
Then we define
$$
f(x)=x^{1/2}+f(a_1)-a_1^{1/2}, \qquad x\in (a_1, a_2]
$$
for some large constant $a_2>a_1$. The constant $a_2$ is chosen so that for large
values of $x$ in $(a_1, a_2]$, the function $f$ behaves like $x^{1/2}$, that is
$f(\lambda x)/f(x)$ is close to $\lambda^{1/2}$ uniformly for $\lambda\in [1, 3]$.
Then we define
$$
f(x)=x^{1/3}+f(a_2)-a_2^{1/3}, \qquad x\in (a_2, a_3]
$$
for some large constant $a_3>a_2$. The constant $a_3$ is chosen so that for large
values of $x$ in $(a_2, a_3]$, the function $f$ behaves like $x^{1/3}$, that is
$f(\lambda x)/f(x)$ is close to
$\lambda^{1/3}$
 uniformly for $\lambda\in [1, 4]$.
We repeat this procedure to define this function inductively.

The function $f$ is an increasing function in OR with upper Matuszewska
index $1/2$ and lower Matuszewska index $1/3$.

Let $\sigma$ be the measure with distribution function
$f$.
Since $\int_{(0, \infty)}(1+t)^{-1}\sigma(dt)<\infty$, $\sigma$ is a
Stieltjes measure. Let
$$
g(\lambda):=\int_{(0, \infty)}\frac1{\lambda + t}\, \sigma(dt)
$$
be the corresponding Stieltjes function. It follows from integration by parts that
$$
g(\lambda)=\int^\infty_0\frac{f(\xi)}{(\lambda+\xi)^2}d\xi.
$$
Using our construction of $f$ we know that
$$
\int^\infty_2\frac{\xi^{1/3}}{(\lambda+\xi)^2}d\xi\le
\int^\infty_2\frac{f(\xi)}{(\lambda+\xi)^2}d\xi\le
\int^\infty_2\frac{\xi^{1/2}}{(\lambda+\xi)^2}d\xi.
$$
Thus it follows from \cite[Lemma 6.3]{WYY} that
$$
c_1\lambda^{-2/3}\le g(\lambda)\le c_2 \lambda^{-1/2},
\quad \lambda\ge 2
$$
for some positive constants $c_1<c_2$.

Modifying the argument of the proof of the de Haan--Stadtm\"uller theorem
(\cite[Theorem 2.10.2]{BGT}) one can show that $g$ is in $OR$ with upper Matuszewska
index $-2/3$ and lower Matuszewska index $-1/2$.
It follows from \cite[Theorem 7.3]{SSV} that $\phi(x):=1/f(x)$ is a
complete Bernstein function.
Thus $\phi$ is a complete Bernstein function in $OR$ with upper Matuszewska
index $1/2$ and lower Matuszewska index $1/3$.

It follows from \cite[Lemma 6.3]{WYY}
that $g$ cannot be comparable with any
regularly varying function at infinity, and therefore $\phi$
cannot be comparable with any
regularly varying function at infinity.

\medskip

{\bf Acknowledgements.} We thank Nick Bingham and Charles Goldie for their
help related to the construction of the function in the last section.

\vspace{.1in}
\begin{singlespace}

\end{singlespace}
\end{doublespace}

\vskip 0.1truein

{\bf Panki Kim}

Department of Mathematical Sciences and Research Institute of
Mathematics,

Seoul National University, San56-1 Shinrim-dong Kwanak-gu, Seoul
151-747, Republic of Korea

E-mail: \texttt{pkim@snu.ac.kr}

\bigskip

\bigskip

{\bf Renming Song}

Department of Mathematics, University of Illinois, Urbana, IL 61801,
USA

E-mail: \texttt{rsong@math.uiuc.edu}

\bigskip

{\bf Zoran Vondra{\v{c}}ek}

Department of Mathematics, University of Zagreb, Zagreb, Croatia

Email: \texttt{vondra@math.hr}

\end{document}